\documentclass{amsart}
\usepackage[all]{xy}
\usepackage{amssymb}
\usepackage{amsmath}
\usepackage{amsmath}
\swapnumbers
\newtheorem{thm}[subsection]{Theorem}
\newtheorem{prop}[subsection]{Proposition}

\newtheorem{cor}[subsection]{Corollary}

\def\d{\succ}
\def\g{\prec}

\def\oo{\omega}

\def\DD{\Delta}

\def\t{\otimes}
\def\pt{\textrm{-}}

\newenvironment{proo}{\begin{trivlist} \item{\emph{Proof.}}}
  {\hfill $\square$ \end{trivlist}}

\newcommand{\Vt}[1]{V^{\otimes #1}}

\def\HH{{\mathcal{H}}}

\def\Prim{\mathrm{Prim\, }}

\def\HHo{\overline {\HH}}
\def\To{\overline {T}}

\def\qa{quasi-shuffle algebra }

\def\ctd{commutative tridendriform }

\begin{document}

\author[J.-L. Loday]{Jean-Louis Loday}
\address{Institut de Recherche Math\'ematique Avanc\'ee\\
    CNRS et Universit\'e Louis Pasteur\\
    7 rue R. Descartes\\
    67084 Strasbourg Cedex, France}
\email{loday@math.u-strasbg.fr}
%\urladdr{www-irma.u-strasbg.fr/{$\sim$}loday/}

\title{On the algebra of quasi-shuffles}
%\alttitle{} 
\subjclass[2000]{16A24, 16W30, 17A30, 18D50, 81R60.}
\keywords{ Bialgebra, generalized bialgebra, Hopf algebra, Cartier-Milnor-Moore, 
Poincar\'e-Birkhoff-Witt, shuffle, quasi-shuffle, stuffle,  dendriform, operad, Fubini number}

%\date{\today}

\begin{abstract} For any commutative algebra $R$ the shuffle product on the tensor module $T(R)$ can be deformed to a new product. It is called the quasi-shuffle algebra, or stuffle algebra, and denoted $T^q(R)$. We show that if $R$ is the polynomial algebra, then $T^q(R)$ is free for some algebraic structure called Commutative TriDendriform (CTD-algebras). This result is part of a structure theorem for CTD-bialgebras which are associative as coalgebras and whose primitive part is commutative. In other words, there is a good triple of operads $(As, CTD, Com)$ analogous to $(Com, As, Lie)$.

In the last part we give a similar interpretation of the quasi-shuffle algebra in  the noncommutative setting.
\end{abstract}

\maketitle

\section*{Introduction} \label{S:int}
The shuffle algebra is a commutative algebra structure on the tensor module $T(V)$ of the vector space $V$. For instance $u \sqcup\!\!\sqcup v = u\t v + v\t u$. If $R$ is a commutative algebra, then the shuffle product on $T(R)$ can de deformed into a new product $*$ called the quasi-shuffle product. For instance
$$ a* b = a\cdot b\ + a\t b + b\t a\ .$$
Observe that this product is not a graded product, but only a filtered product, since $a\t b\ + b\t a \in R^{\t 2}$ and $a\cdot b \in R$.

The quasi-shuffle algebra appears in different places in mathematics, for instance in the Multi-Zeta Values topic, cf.~for instance \cite{Ca2, H, IKZ}, where it plays an important role.

The aim of this paper is to show, first, that the quasi-shuffle algebra over the polynomial algebra has a universal property: it is the free ``Commutative TriDendriform algebra''. This notion of CTD-algebra is a commutative version of the notion of dendriform algebras, which plays a role in several themes (cf.~\cite{BFK, LR1, V}). A CTD-algebra is determined by two operations, denoted $\g$ (left) and $\cdot$ (dot), related by the axioms:
$$ x\cdot y = y\cdot x$$
and
\begin{eqnarray*}
(x\g y)\g z &=& x\g (y\g z + z\g y + y\cdot z),\\
(x\cdot y)\g z &=& x\cdot (y\g z ),\\
(x\cdot y)\cdot z &=& x\cdot (y\cdot z ).
\end{eqnarray*}
We show that the functor $T^q : Com\textrm{-alg} \to CTD\textrm{-alg} $ is left adjoint to the forgetful functor (which forgets $\g$). 
A quasi-shuffle algebra is not only a CTD-algebra, but it is even a CTD-bialgebra. In the second part we state and prove a structure theorem for connected CTD-bialgebras, cf.~\ref{asctdcom}. It means that the triple of operads
$$(As, CTD, Com)$$
is good in the sense that it satisfies the analogues of the Poincar\'e-Birkhoff-Witt theorem and the Cartier-Milnor-Moore theorem satisfied by the triple of operads $(Com, As, Lie)$.

In the last part we state and prove a non-commutative analogue of the universal property: the quasi-shuffle algebra over the tensor algebra is free in the category of involutive tridendriform algebras.\\

\noindent {\bf Notation.} In this paper $K$ is a field and all vector spaces 
are over $K$. Its unit is denoted by 1. 
%The vector space spanned by the elements of a set $X$ is denoted $K[X]$. 
The tensor product  of vector spaces over $K$
is denoted by $\t$. The tensor product of $n$ copies of the space $V$ is
denoted $\Vt n$. For $v_i\in V$ the element $v_1\t \cdots \t v_n$ of 
$\Vt n$ is denoted by the concatenation of the elements: $v_1 \ldots v_n$.
%A linear map $\Vt n \to V$ is called an {\it $n$-ary operation} on $V$ 
%and a linear map $ V\to \Vt n$ is called an {\it
%$n$-ary cooperation} on $V$.
The \emph{tensor module} over $V$ is the direct sum
$$T(V) := K 1\oplus V\oplus V^{\t 2}\oplus \cdots \oplus V^{\t n} \oplus \cdots $$
and the \emph{reduced tensor module} is $\To(V) := T(V)/K1$.
%$$\To(V) := V\oplus V^{\t 2}\oplus \cdots \oplus V^{\t n} \oplus \cdots .$$
The \emph{symmetric module} over $V$ is the direct sum
$$S(V) := K 1\oplus V\oplus S^2(V) \oplus \cdots \oplus S^n(V) \oplus \cdots $$
and the \emph{reduced symmetric module} is $\bar S(V) := S(V)/K1$,
%$$\bar S(V) := V\oplus S^2(V) \oplus \cdots \oplus S^n(V) \oplus \cdots ,$$
where $S^n(V)$ is the quotient of $V^{\t n}$ by the action of the symmetric group. We still denote by $v_1 \ldots v_n$ the image in $S^n(V)$ of $v_1 \ldots v_n\in V^{\t n}$.
If $V$ is generated by $x_1,\ldots, x_n$, then $S(V)$ (resp. $T(V)$) can be identified with the polynomials (resp. noncommutative polynomials) in $n$ variables.

\noindent{\bf Acknowledgement} I thank Muriel Livernet for comments on the first version of this paper. This work has been partially supported by the ``Agence Nationale pour la Recherche''.

 \section{Quasi-shuffle algebra}\label{S:quasishuffle}

\subsection{Definition, notation}\label{quasishuffle} For any associative (not necessarily unital) algebra $(R, \cdot)$ the \emph{quasi-shuffle algebra}, or \emph{stuffle algebra}, on $R$ is the tensor module $T(R)$ equipped with the associative unital product $*$ uniquely determined by the relation:
$$ax * by = (a\cdot b)( x * y ) +a( x * by )  +  b(a x * y ) $$
for any $a,b\in R$, any $x, y  \in T(R)$, 1 being a unit for $*$.

The product $*$ is associative and unital, we denote by $T^q(R):= (T(R), *)$ the \qa on $R$. It is immediate to check that if $(R, \cdot)$ is commutative, then so is $T^q(R)$. 

Explicitly the product $*$ is the following. Let $\gamma$ be a piece-wise linear path in the plane from $(0,0)$ to $(p,q)$ using only the steps $(0,1), (1,0)$ and $(1,1)$. Define $\gamma(a_{1}\ldots a_{p}, a_{p+1}\ldots a_{p+q})$ to be the element of $R^{\t r}, r\leq p+q$, obtained by writing $a_{1}$ if the first step is $(0,1)$, $a_{p+1}$ if the first step is $(1,0)$, and $a_{1}\cdot a_{p+1}$ if the first step is $(1,1)$, etc. until the last step. Observe that, if there is no $(1,1)$ step, then we obtain a shuffle. If the path is the diagonal ($p=q$), then we get $a_{1}\cdot a_{p+1}\cdot \ldots \cdot a_{p}\cdot a_{2p}\in R$. With this notation we have
$$a_{1}\ldots a_{p} * a_{p+1}\ldots a_{p+q}=\sum_{\gamma}\gamma(a_{1}\ldots a_{p}, a_{p+1}\ldots a_{p+q}).$$

The elements in the sum are sometimes called ``mixable shuffles'', cf.~\cite{GK}.

\subsection{Examples} 1. If the product $\cdot$ is zero (recall that $R$ need not be unital), then $T^q(R)$ is the shuffle algebra on the vector space $R$.

\noindent 2. Let $y_i$, $i\geq 1$, be a set of variables and let $R_0:=K\langle y_1, y_2, \ldots \rangle$ be the noncommutative polynomial algebra on a countable number of variables. Define a product $*$ on $R_0$ by 
$$y_k\oo  * y_{k'}\oo '  = y_{k+k'}( \oo  * \oo ' ) +y_k( \oo  * y_{k'}\oo ' )  +  y_{k'}(y_k \oo  * \oo ' ) ,$$
see formula (94) of \cite{Ca2}. Then it is easy to check that $R_0$ is the \qa on the algebra of non-constant polynomials in one variable. 

Let us denote by $T^c(V)$ the cofree coalgebra on the vector space $V$. As a vector space it is the tensor module $T(V)$. The comultiplication is  the \emph{deconcatenation} given by 
$$\DD(v_1\ldots v_n) := \sum_{i=0}^n v_1\ldots v_i \t v_{i+1}\ldots v_n.$$

\begin{prop} The unique coalgebra morphism
$$T^c(R) \t T^c(R) \to T^c(R)$$
induced by the algebra morphism $\cdot : R\t R \to R$ is the quasi-shuffle product $*$.
\end{prop}
\begin{proo} Since $T^c(R) $ is cofree in the category of connected coassociative coalgebras, there is a unique coalgebra morphism $\mu :T^c(R) \t T^c(R) \to T^c(R)$ which makes the following diagram commutative:
$${\xymatrix{
T^c(R)\t T^c(R)\ar[r]^-{\mu}\ar@{->>}[d] & T^c(R)\ar@{->>}[d]  \\
R\t R\ar[r]^{\cdot}& R\\
}}$$
Since $*$ is a coalgebra morphism, in order to show that $\mu=*$, it is sufficient to show that the two composites $T^c(R)\t T^c(R)\to R\t R \to R$ and $T^c(R)\t T^c(R)\to T^c(R) \to R$ are equal. By definition $ a* b = ab + ba + a\cdot b$ for any $a,b\in R$, hence the projection of this element in $R$ is $a\cdot b$ and we are done.
\end{proo}

\subsection{Remark} If, instead of starting with an algebra $R$, we start with a coalgebra $C$, then the dual construction gives rise to a bialgebra structure on the tensor algebra $T(C)$. This construction is used in \cite{BS} to construct renormalization Hopf algebras in various settings. See also \cite{Ca1, Man} where such objects were studied.

\subsection{Operations on the \qa } Let us denote by $\g$ (called left), by $\d$ (called right) and by $\cdot$ (called dot) the three  binary operations on $T(R)$ uniquely determined by the following requirements (inductive definition):

$$\begin{array}{rclccc}
ax \cdot by  &=& (a\cdot b)( x * y ) & and & 1\cdot x=0& x\cdot 1 =0,\\
ax \g by  &=& a( x * by ) & and & 1\g x=0& x\g1 =x,\\
ax \d by  &=&  b(a x * y ) & and & 1\d x=x& x\d1 =0,
\end{array}$$
where, for any $x,y\in T(R)$, we have defined
$$x*y :=   x\g y +  x\d y + x\cdot y .$$
It is immediate to check that this operation is the same as the quasi-shuffle introduced in \ref{quasishuffle}.
\begin{prop}\label{threeop}  The three operations  $\g$, $\d$ and $\cdot$ on $T(R)$ satisfy the following relations:
\begin{eqnarray}
(x \g y) \g z &=& x \g (y * z)\\
(x \d y) \g z &=& x \d (y \g z)\\
(x * y) \d z &=& x \d (y \d z)\\
(x \cdot y) \g z &=& x \cdot (y \g z)\\
(x \g y) \cdot z &=& x \cdot (y \d z)\\
(x \d y) \cdot z &=& x \d (y \cdot z)\\
(x \cdot y) \cdot z &=& x \cdot (y \cdot z).
\end{eqnarray}
\end{prop}

\begin{proo} The proof is by induction on the degree of the elements (an element of $V^{\t n}$ is of degree $n$). First, we prove the seven relations when the elements are in $R$.\\

\noindent
1. $(a\g b)\g c = (ab)\g c = a(b*c)= a\g (b*c)$,\\
2. $(a\d b)\g c = (ba)\g c = b(a*c) = a\d (bc) = a\d (b\g c)$,\\
3. $(a*b)\d c = c(a*b) = a\d (cb) = a\d (b \d c)$,\\
4. $(a\cdot b) \g c = (a\cdot b)c = a\cdot (bc)= a\cdot (b\g c)$,\\
5. $(a\g b)\cdot c = (ab)\cdot c = (a\cdot c)b =a\cdot (cb) = a\cdot (b\d c)$,\\
6. $(a\d b) \cdot c = (ba) \cdot  c = (b\cdot c)a = a\d (b\cdot c)$,\\
7. $(a\cdot b)\cdot c = a\cdot (b\cdot c)$.\\

Second, we show by induction the first relation for elements in $T(R)$. On the one hand we have
\begin{eqnarray*}
(ax \g y) \g z &=& a(x*y)\g z\\
& = & a((x*y)*z)
\end{eqnarray*}
 On the other hand we have
\begin{eqnarray*}
 ax \g (y * z) &=&a(x*(y*z)).
 \end{eqnarray*}
Since $*$ is associative by the inductive assumption, we get the first relation $(x \g y) \g z = x \g (y * z)$. The other relations are proved similarly.
\end{proo}

\begin{prop}\label{threerel} If $(R, \cdot)$ is commutative, then  $x\cdot y = y \cdot x$, $x\d y=y\g x$ and  $x*y=y*x$ for any $x,y\in T(R)$. Under this hypothesis the seven relations above come down to the following three:
\begin{eqnarray*}
(x\g y)\g z &=& x\g (y* z),\\
(x \cdot y) \g z &=& x \cdot (y \g z),\\
(x\cdot y)\cdot z &=& x\cdot (y\cdot z ).
\end{eqnarray*}
\end{prop}
\begin{proo} The first part of the statement is immediate. For the second part: by the symmetry hypothesis the 7 relations can be written in terms of $\g$ and $\cdot$. Then relation (3) is equivalent to relation (1), and relation (2) is implied by relation (1). Relation (6) is equivalent to relation (4) and relation (5) is implied by relation (4). So we are left with relations (1), (4) and (7) which can be written as in the statement because of the symmetry hypothesis.
\end{proo}

\subsection{Remark} A ${\bf B}_{\infty}$-algebra $R$ is a vector space equipped with $(p+q)$-ary operations $M_{pq}$ verifying the relations which enable us to put a bialgebra structure on $T^c(R)$. If $M_{pq}=0$ for $p\geq 2$, then it is called a \emph{brace algebra}. If $M_{pq}=0$ for $p\geq 2$ and  $q\geq 2$, then only $M_{11}$ survives and $R$ is an associative algebra. The general procedure to construct the bialgebra out of the ${\bf B}_{\infty}$-algebra (cf.~for instance \cite{LR2}), gives precisely the quasi-shuffle algebra in the latter case.

\section{Commutative TriDendriform algebra}

The properties of the \qa shown in the previous section suggest the notion of \ctd algebra that we now introduce. We show that the free CTD-algebra is the quasi-shuffle algebra over the polynomials.

The notion of dendriform algebra was first introduce in \cite{Lod1} in connection with some problems in algebraic $K$-theory, and the notion of tridendriform algebra was introduced in \cite{LR1}.

\subsection{Definition} A \emph{\ctd algebra} or $ComTriDend$-algebra (CTD-algebra for short)  is a vector space $A$ equipped with two operations $\g$ and $\cdot$ satisfying:
$$x\cdot y = y\cdot x$$
and
\begin{eqnarray}
(x\g y)\g z &=& x\g (y\g z + z\g y + y\cdot z),\\
(x\cdot y)\g z &=& x\cdot (y\g z ),\\
(x\cdot y)\cdot z &=& x\cdot (y\cdot z ).
\end{eqnarray}

It is immediate to check from these relations that the product $*$, defined as $x*y := x\g y + y\g x  + x\cdot y$, is associative and commutative.

It will prove helpful to introduce the notion of \emph{ unital CTD-algebra} which is a vector space of the form $A = K1\oplus \bar A$ where $\bar A$ is CTD-algebra and the operations $\g$ and $\cdot$ are extended to $1$ by the following formulas:
$$  1\cdot x=0=x\cdot 1, 1\g x=0, x\g1 =x \textrm{ for any } x\in \bar A.$$
Observe that $1\g 1$ and $1\cdot 1$ are not defined but $1$ is a unit for $*$. Therefore $\bar A$ is the augmentation ideal of $A$.

\subsection{Examples}\label{examp} (a) From Proposition \ref{threerel} it is clear that the \qa over a commutative ring is a CTD-algebra. In particular the algebra of shuffles and  $T^q(\bar S(V))$ are CTD-algebras.

(b) Let $(R, \cdot)$ be a commutative algebra and let $P:R\to R$ be a \emph{Rota-Baxter operator}, i.e.~a linear map which satisfies the identity:
$$P(a)\cdot P(b) = P(a\cdot P(b) + P(a)\cdot  b + a\cdot b).$$
If we define $a\g b := a\cdot P(b)$, then one can check easily that $(R, \g, \cdot)$ is a $CTD$-algebra, cf.~\cite{A, Ca1, EFG}. Moreover the operator $P$ becomes an associative algebra morphism: $P: (R,*) \to (R, \cdot)$.

\subsection{Free CTD-algebra} Let $V$ be a vector space. By definition the \emph{free CTD-algebra} over $V$ is a CTD-algebra $CTD(V)$ equipped with a linear map $\iota: V \to CTD(V)$ which satisfies the following universal condition:

any linear map $\phi :V\to A$, where $A$ is a CTD-algebra, admits a unique extension $\Phi : CTD(V) \to A$ as morphism of CTD-algebras:

$${\xymatrix{
CTD(V)\ar[r]^-{\Phi} & A \\
V\ar[u]^{\iota}\ar[ur]_{\phi}& \\
}}$$
In the unital case it is assumed that the image of $\phi$ lies in the augmentation ideal of $A$. 

\subsection{The operad CTD}  Let us work over a field of characteristic zero field. 
 Since the relations of a CTD-algebra are multilinear, the $n$-multilinear part $CTD(n)$ of $CTD(Kx_{1}\oplus \cdots \oplus Kx_{n})$ inherits a structure of $S_{n}$-module from the action of the symmetric group on the set of variables $\{x_{1},\ldots ,x_{n}\}$, and we have
$$CTD(V) = \bigoplus_{n}CTD(n)\t _{S_{n}}V^{\t n}.$$
 So the functor $CTD$ is a Schur functor. It is called the operad of the CTD-algebras (cf.~the first Chapter of \cite{Lod3} for instance).
As a representation of $S_{n}$, $CTD(n)$ is spanned by the ordered-unordered partitions of  $\{1,\ldots ,n\}$.
The dimension $d_{n}$ of $CTD(n)$ is sometimes called the \emph{Fubini number} (sequence number A000670 in the On-Line Encyclopedia of Integer Sequences):
 $$\begin{array}{cccccccc}
n& 1 & 2 & 3 & 4 & 5 & 6 &  \cdots \\
d_{n} & 1 & 3 & 13 & 75 & 541 & 4683 &\cdots
\end{array}
$$
It can be shown combinatorially that the generating series $f^{CTD}(x):=\sum_{n\geq 1} \frac{d_{n}}{n!} x^n$ is given by:
$$f^{CTD}(x)= \frac{\exp (x)-1}{2- \exp(x)} .$$
It will be a consequence of our main result, cf.~\ref{genser}.

\begin{thm}\label{Phi} The free unital CTD-algebra on the vector space $V$ is isomorphic to the \qa over the reduced symmetric algebra on $V$ (non-constant polynomials):
$$\Phi :CTD(V)\cong T^q(\bar S(V)).$$
\end{thm}
\begin{proo} Since by Proposition \ref{threerel} the \qa $T^q(\bar S(V))$ is a CTD-algebra and since it contains $V$, the inclusion map $\iota$ induces a CTD-morphism $\Phi(V) : CTD(V)\to T^q(\bar S(V)).$ It is clear from the definition of the operations $\g$ and $\cdot$ in $T^q(\bar S(V))$ that $\Phi$ is surjective. So we need only to check that it is injective. From the naturality of $\Phi$ it suffices to check injectivity on the multilinear parts (since both functors are Schur functors). We denote them by $CTD(n)$ and $T^q\bar S (n)$ respectively.

Let $V$ be an $n$-dimensional vector space, spanned by $v_{1}, \ldots ,v_{n}$.  In $T^q(\bar S(V))$ the multilinear part of degree $n$ admits as a basis the vectors $$(w_{1}\cdots w_{i_{1}})
(w_{i_{1}+1}\cdots w_{i_{1}+i_{2}})\cdots ( \cdots\  \cdots w_{i_{1}+ \cdots + i_{k}})$$
 where the $w_{i}$'s are a permutation of the $v_{i}$'s and the subwords in parenthesis are in increasing order. We call such an element an \emph{ordered-unordered partition}. Example: $(v_{2}v_{5})(v_{1}v_{4})(v_{3})$. It suffices to show that $\dim CTD(n) \leq \dim T^q\bar S (n)=:d_{n}$. 

First, we show that any monomial in $CTD(V)$ can be written as
$$x_{1}\g (x_{2}\g (\cdots (x_{k-1}\g x_{k})\cdots))$$
where $x_{i}$ is a monomial of the form $v_{i_{1}}\cdot\ldots \cdot v_{i_{l}}$ where the indices $i_{j}$ are in increasing order.

Any monomial of degree larger than 1 is of the form $X\g Y$ or $X\cdot Y$ for elements of strictly smaller degree. We freely use the fact that the operation $\cdot$ is associative and commutative. We work by induction on the degree of the components.

\noindent (a) Case $X\g Y$. 

(a1) If $X=v_{i_{1}}\cdot\ldots \cdot v_{i_{l}}$, we are done since the degree of $Y$ is strictly smaller than the degree of $X\g Y$.

(a2) Otherwise, either $X=X_{1} \cdot X_{2}$ where the symbol $\g$ appears in $X_{1}$ or $X_{2}$, or $X=X_{1} \g X_{2}$. In the second case we apply relation (8) and we are back to case (a) with the degree of the second component being strictly smaller. In the first case, on can suppose that $X_{1}= X'_{1} \g X''_{1}$. We apply relation (9) to obtain $((X'_{1} \cdot X_{2})\g   X''_{1}) \g Y$ to which we apply relation (8). Again, we are back to case (a) with the degree of the first component strictly smaller.

\noindent (b) Case $X\cdot Y$. If the product $\g$ does not appear in $X$ nor $Y$ we are done. Otherwise, assume that $X=X_{1} \g X_{2}$. We use relation (9) to get back to case (a) as before and we are done.

We have shown that $CTD(n)$ is spanned by $d_{n}$ elements, where $d_{n}$ is the number of ordered-unordered partitions. Since $CTD(n) \to T^q\bar S (n)$ is surjective we are done.
\end{proo}

\subsection{Universal enveloping CTD-algebra} Let $A$ be a CTD-algebra. Ignoring the left operation yields a forgetful functor from the category of CTD-algebras to the category of commutative algebras (not necessarily unital). We denote by $U_{CTD}$ the left adjoint of this forgetful functor:
$$U_{CTD} : Com\textrm {-alg}\to CTD\textrm {-alg}.$$
Explicitly, for any commutative algebra $R$,  the CTD-algebra $U_{CTD}(R)$ is given by
$$U_{CTD}(R)= CTD(R)/\sim$$
where the equivalence relation $\sim$ consists in identifying,  for any $a,b\in R$,  the product in $R$, $ab\in R\subset CTD(R)$, with the dot-product in $CTD(R)$, $a\cdot b\in CTD(R)$. Observe that only the vector space structure of $R$ is used to construct the free CTD-algebra $CTD(R)$. The commutative structure of $R$ is used when making the quotient.

\begin{thm}\label{thm:uctd} For any (not necessarily unital) commutative algebra $R$ there is an isomorphism of CTD-algebras:
$$U_{CTD}(R) \cong T^q(R).$$
\end{thm}
\begin{proo} The inclusion $(R,\cdot) \to (T^q(R),\cdot)$ determines an adjoint map $U_{CTD}(R) \to T^q(R)$ which is a CTD-morphism. In order to prove that it is an isomorphism it suffices, by a classical argument, cf.\cite{Qui} Appendix B, to prove it for the free commutative algebra $R= \bar S(V)$. Since $U_{CTD}$ and $\bar S$ are both left adjoint functors, their composite is left adjoint to the forgetful functor $CTD\textrm{-alg}\to Vect$ and therefore we have $U_{CTD}(\bar S(V)) = CTD(V)$. The restriction of the adjoint map $CTD(V)=U_{CTD}(\bar S(V)) \to T^q(\bar S(V))$ to $V$ is simply the inclusion $V\to \bar S(V) \to  T^q(\bar S(V))$, therefore this adjoint map is the isomorphism $\Phi$ of Theorem \ref{Phi}.
\end{proo}

\section {CTD-bialgebras} We use the technique developed in \cite{Lod2} to show that the free CTD-algebra is in fact a Hopf algebra. We introduce the notion of CTD-bialgebra, which is a CTD-algebra equipped with a compatible coproduct. Quasi-shuflle algebras are CTD-bialgebras.
In this section we always assume unitality.

\subsection{Tensor product of CTD-algebras}\label{tensorproduct} Let $A$ and $B$ be two CTD-algebras. One can put a CTD-algebra structure on $A\t B$ as follows:
\begin{eqnarray*}
(a\t b) \g (a'\t b') & = & (a\g a') \t (b*b'), \\
(a\t b) \cdot (a'\t b') & = & (a\cdot a') \t (b*b') .
\end{eqnarray*}
Since the elements $1\g 1$ and $1\cdot 1$ do not make sense, we take the following convention: 
\begin{eqnarray*}
1\g 1  \t b* b' & := & 1 \t (b\g b') ,\\
1\cdot 1  \t b* b'  & := & 1 \t (b\cdot b') .
\end{eqnarray*}
This trick is due to M. Ronco. Observe that if $A, B$ and $C$ are three CTD-algebras, then the two CTD-algebra structures that one can put on $A\t B \t C$ are the same.

\subsection{Bialgebra structure on the free CTD-algebra}\label{bialgebraCTD} Let CTD(V) be the free (unital) CTD-algebra on $V$. Since $CTD(V)\t CTD(V)$ is a CTD-algebra there is a unique CTD-morphism
$$\Delta: CTD(V) \to CTD(V)\t CTD(V)$$
which extends the linear map $V \to CTD(V)\t CTD(V)$ given by $v\mapsto v\t 1 + 1 \t v$.

As in \cite{Lod2} we see that $\DD$ is a coassociative counital morphism.

This example justifies the following definition.

\subsection {$As^c\pt CTD$-bialgebra}\label{ctdbialgebra} By definition an \emph{$As^c\pt CTD$-bialgebra}, or CTD-bialgebra for short, is a vector space $\HH = K1\oplus \HHo$ equipped with a structure of unital CTD-algebra, a structure of coassociative counital coalgebra, and these two structures are related by the following \emph{compatibility relation}:
\begin{eqnarray*}
\DD(x\g y) &= & x_{(1)}\g y_{(1)}\t  x_{(2)}* y_{(2)},\\
\DD(x\cdot y)& = & x_{(1)}\cdot y_{(1)}\t  x_{(2)}* y_{(2)},
\end{eqnarray*}
where we use the notation $\DD(x) = x_{(1)}\t x_{(2)}$. Here we adopt the convention stated in \ref{tensorproduct} and $\DD(1) = 1\t 1$.

Following Quillen (cf.~\cite{Qui}, p.~282) we say that a coaugmented coalgebra $\HH$ is 
{\it connected} if $\HH = \bigcup _{r\geq 0} F_r\HH$
where $F_r\HH$ is the coradical filtration of $\HH$ defined recursively by the 
formulas
\begin{eqnarray*}
F_0\HH &:=& K 1, \\
F_r\HH & :=& \{ x\in \HH \mid  \DD  (x)- x\t 1 - 1 \t x \in F_{r-1}\HH \t F_{r-1}\HH \}\ .
\end{eqnarray*}

From the fact that $CTD(V)$ is an $As^c\pt CTD$-bialgebra, cf. \ref{bialgebraCTD}, it follows that 
$U_{CTD}(R)$ is also an $As^c\pt CTD$-bialgebra.

\begin{prop}\label{prim} For any CTD-algebra $\HH$ the dot operation is stable on the primitive part $\Prim \HH :=\{x\in \HH\ \vert \ \DD(x) - x\t 1 - 1\t x=0\}$. So $\Prim \HH$ is a commutative algebra.
\end{prop}
\begin{proo} It suffices to compute for $x$ and $y$ primitive:
\begin{eqnarray*}
\DD(x\cdot y) &=& (x\t 1 + 1 \t x) \cdot (y\t 1 + 1 \t y) \\
& = & x\cdot y \t 1*1 +  x\cdot 1 \t 1*y +  1\cdot y \t x*1 +  1\cdot 1 \t x*y\\
& = & x\cdot y \t 1 + 0 + 0 +  1 \t x\cdot y
\end{eqnarray*}
\end{proo}
\begin{prop}\label{isophi} The isomorphism $\Phi: U_{CTD}(R) \to T^q(R)$, cf.~Theorem \ref{Phi}, is compatible with the coproduct, and so is an isomorphism of CTD-bialgebras.
\end{prop}
\begin{proo} Since the coproduct on $U_{CTD}(R)$ is induced by the coproduct on $CTD(R)$, it suffices to show that $\Phi: CTD(R) \to T^q(R)$ is compatible with the coproduct.  From the definition of the coproduct on $CTD(R)$, cf.~\ref{bialgebraCTD}, it suffices to show that the coproduct on $T^q(R)$, that is deconcatenation, satisfies the compatibility relations of \ref{ctdbialgebra}. We adopt the notation $a(x\t y) := (ax)\t y$ for $a\in R$ and $x\t y\in T(R)\t T(R)$ and we remark that $\DD(az) = a\DD(z) + 1\t az$.
For the dot operation we already verified the compatibility for elements of $R$ in Proposition \ref{prim}. For elements of $T^q(R)$ we have on the one hand:
\begin{eqnarray*}
\DD(ax\cdot by) & = & \DD((a\cdot b)(x*y))\\
&=& (a\cdot b)\DD(x*y)+ 1\t (a\cdot b)(x*y)
\end{eqnarray*}
On the other hand we get:
\begin{eqnarray*}
\DD(ax)\cdot \DD(by) & = &(a\DD(x)+ 1\t ax)\cdot (b\DD(y)\t 1 +  1 \t by)\\
& = &(a\cdot b)\DD(x)*\DD(y)+  1 \t ax\cdot by\\
&=& (a\cdot b)\DD(x*y)+ 1\t (a\cdot b)(x*y)
\end{eqnarray*}
Therefore we have proved that $\DD(ax\cdot by)= \DD(ax)\cdot \DD(by) $ as expected.

For the left operation, we first verify it on the elements of $R$.

One hand hand we get:
\begin{eqnarray*}
\DD(a\g b) & = & \DD(ab)\\
& = & ab\t 1 + a\t b + 1 \t ab.
\end{eqnarray*}
On the other hand we get:
\begin{eqnarray*}
\DD(a)\g \DD(b) & = &(a\t 1 +  1 \t a)\g (b\t 1 +  1 \t b)\\
&=& a\g b \t 1 + a\g 1 \t 1*b + 1\g b \t a*1 + 1 \g 1 \t a*b\\
&=& a b \t 1 + a \t 1*b + 0+ 1\t a\g b\\
& = & ab\t 1 + a\t b + 1 \t ab.
\end{eqnarray*}
We now prove this relation for any elements of $T^q(R)$ by induction. 
On the one hand we get:
\begin{eqnarray*}
\DD(ax\g by) & = & \DD(a(x*by))\\
&=& a\DD(x*by)+ 1\t a(x*by)\\
&=& a(\DD(x)*\DD(by))+ 1\t a(x*by)\\
&=& a(\DD(x)*b\DD(y)) +a(\DD(x)*(1\t by)) + 1\t a(x*by)\\
\end{eqnarray*}
On the other hand we get:
\begin{eqnarray*}
\DD(ax)\g \DD(by) & = &(a\DD(x)+ 1\t ax)\g (b\DD(y) +  1 \t by)\\
& = &a\DD(x)\g b\DD(y)+ a\DD(x)\g ( 1 \t by)+ 1\t ax\g by\\
& = &a(\DD(x)* b\DD(y))+ a(\DD(x)*( 1 \t by))+ 1\t a(x*by).
\end{eqnarray*}
Therefore we have proved that $\DD(ax\g by)= \DD(ax)\g \DD(by) $ as expected.
\end{proo}

\section{Structure theorem for CTD-bialgebras}
We prove a structure theorem for $As^c\pt CTD$-bialgebras analogous to the structure theorem for cocommutative Hopf algebras (PBW+CMM theorem), cf.~\cite{MM, Qui}. Here the structure of the primitive part is the commutative structure (instead of the Lie structure in the classical case). Therefore $(As, CTD, Com)$ is a good triple of operads in the sense of \cite{Lod3}.

\begin{thm}\label{asctdcom} We suppose that the ground field $K$ is of characteristic zero. If $\HH$ is a CTD-bialgebra over $K$, then 
the following are equivalent:\\
\noindent (a) $\HH$ is a connected CTD-bialgebra,\\
\noindent (b) $\HH$ is isomorphic to $U_{CTD}({\Prim} {\HH})$ as a CTD-bialgebra,\\
\noindent (c) $\HH$ is cofree among the connected coalgebras: $\HH\cong T^c({\Prim} {\HH})$.
\end{thm}
\begin{proo} We apply the main theorem of \cite{Lod3} to the $As^c\pt CTD$-bialgebra type. In order to get the equivalence of the three assertions, it is sufficient to verify three hypotheses, which, in our case, read as follows:

\noindent (H0) the compatibility relation is distributive,

\noindent (H1) the free CTD-algebra is naturally a  $As^c\pt CTD$-bialgebra,

\noindent (H2epi) the natural coalgebra map $\varphi(V): CTD(V) \to T^c(V)$ admits a natural splitting.

Hypothesis (H0) is immediate by direct inspection, cf.~\ref{ctdbialgebra}. Hypothesis (H1) is immediate by construction, cf.~\ref{bialgebraCTD}. Before proving hypothesis (H2epi), let us recall the construction of $\varphi(V)$. Since by hypothesis (H1), $CTD(V)$ is a coassociative coalgebra, the natural projection $CTD(V) \to V$ induces a unique coalgebra map $\varphi(V): CTD(V) \to T^c(V)$. By Theorem \ref{Phi} there is an isomorphism 
$CTD(V)\cong T(\bar S(V))$. The projection $CTD(V) \to V$ is the composite of the projections 
$T(\bar S(V))\to \bar S(V)\to V$. The splitting $s: T^c(V) \to T(\bar S(V))$ is simply induced by the inclusion $V \to \bar S(V)$ to which we apply the functor $T^c$.
\end{proo}

\begin{cor}\label{genser} There is an isomorphism of Schur functors $CTD = As \circ Com$ and therefore:
$$f^{CTD}(x)= \frac{\exp (x)-1}{2- \exp(x)} .$$
\end{cor}
\begin{proo} Since the composition of two left adjoint functors is still a left adjoint functor, there is an isomorphism $CTD(V)\cong U_{CTD}(Com(V))$ and therefore an isomorphism $CTD(V)\cong T^c(Com(V))$. As a consequence we get $f^{CTD}= f^{As}\circ f^{Com}$, since the generating series $f^{As}(x)$ of $\To$ is the geometric series $\frac{x}{1-x}$ and  the generating series $f^{Com}(x)$ of $Com$, that is $\overline S$, is the exponential series $\exp(x) -1$. \end{proo}

\subsection{Remark} If we are interested only in the Hopf algebra structure, then $T^q(R)$ is isomorphic to the shuffle algebra $T^{sh}(R)$ of the vector space $R$. Indeed the dual Hopf algebra of $T^q(R)$ is connected and cocommutative. So, by the Cartier-Milnor-Moore theorem, it is isomorphic to the universal enveloping algebra of its primitive part. Since this primitive part is the free Lie algebra $Lie(R)$, we get $U(Lie(R))= T(R)$ where only the vector space structure of $R$ is involved. Since the dual of the tensor Hopf algebra is the shuffle Hopf algebra, we are done.

\subsection{Triple of operads} Theorem \ref{asctdcom} provides a new example of ``good triple of operads" as defined in \cite{Lod3}:
$$(As, CTD, Com).$$
Its quotient triple is $(As, Zinb, Vect)$, where $Zinb$ is the operad of Zinbiel algebras. Recall that a Zinbiel algebra (cf.~\cite{Lod2}) is a vector space equipped with a binary operation $\g$ verifying
$$(x\g y)\g z = x\g (y\g z + z\g y ).$$
Since $(As, CTD, Com)$ is a good triple, so is $(As, Zinb, Vect)$ by a Theorem of \cite{Lod3}. It is also a consequence of a theorem due to M. Ronco which asserts that the triple $(As, Dend, Brace)$ is good, where $Brace$ is the operad of brace algebras, cf.~\cite{R1, R2}. For a self-contained proof concerning $(As, Zinb, Vect)$, see \cite{Bu}.

\subsection{Graded version} Until now we have worked in the symmetric category of vector spaces over $K$. We could as well work in the symmetric category of \emph{graded} vector spaces over $K$, taking the sign into account in the symmetry operator: $x\t y \mapsto (-1)^{\vert x\vert \vert y\vert }y\t x$. In this context the first relation of a graded CTD-algebra reads
$$(x\g y)\g z = x\g (y\g z +  (-1)^{\vert y\vert \vert z\vert }z\g y + y\cdot z)$$
on homogeneous elements. The other two relations are unchanged and the dot operation is supposed to be graded-symmetric.

All the preceding results can be written in the graded framework, provided that the functors are replaced by their graded analogues. For instance, the symmetric functor $S$ has to be replaced by the graded-symmetric functor $\Lambda$ given by
$$\Lambda(V^{ev}\oplus V^{odd}) = S(V^{ev})\t E(V^{odd})$$
where $E$ is the exterior algebra functor.

\section{Structure of non-commutative quasi-shuffle algebras}

\subsection{Definition}\cite{LR1} A \emph{tridendriform algebra} (also called dendriform trialgebra) is a vector space $A$ over $K$ equipped with three binary operations $\g, \d, \cdot$ satisfying the following relations:
\begin{eqnarray*}
(x \g y) \g z &=& x \g (y * z)\\
(x \d y) \g z &=& x \d (y \g z)\\
(x * y) \d z &=& x \d (y \d z)\\
(x \d y) \cdot z &=& x \d (y \cdot z)\\
(x \g y) \cdot z &=& x \cdot (y \d z)\\
(x \cdot y) \g z &=& x \cdot (y \g z)\\
(x \cdot y) \cdot z &=& x \cdot (y \cdot z),\\
\end{eqnarray*}
where $x*y:= x\g y + x \d y + x\cdot y$.
Observe that the product $*$ is associative. As in the commutative case (cf. \cite{Lod2}) we can introduce a unit by requiring:
\begin{eqnarray*}
&1\cdot x=0=x\cdot 1,&\\
 &1\g x=0,\quad x\g1 =x,&\\
 & 1\d x=x,\quad x\d1 =0.&
\end{eqnarray*}

\subsection{Examples} (a) The quasi-shuffle algebra $T^q(R)$ over the (nonunital) associative algebra $R$ is a tridendriform algebra by Proposition \ref{threeop}.

\noindent (b) The free tridendriform algebra on one generator can be described explicitly in terms of planar trees, cf. \cite{LR2}.

\noindent (c) Let $(R,\cdot)$ be a (nonunital) associative algebra and let $P$ be a Baxter operator on $R$, cf.~\ref{examp}. If we put $a\g b:=a\cdot P(b)$ and $a\d b := P(a)\cdot b$, then $(R,\g,\d,\cdot)$ is a tridendriform algebra.

\noindent (d) A CDT-algebra is a particular case of tridendriform algebra as proved in \ref{threerel}.

\subsection{Tridendriform algebra with involution} An \emph{involution} on a tridendriform  algebra $A$ is a linear map $a\mapsto \iota(a)$ verifying the following conditions:
\begin{eqnarray*}
\iota (x\g y)  &=& \iota (y) \d \iota (x)\\
\iota (x\d y)  &=& \iota (y) \g \iota (x)\\
\iota (x\cdot y)  &=& \iota (y) \cdot \iota (x)
\end{eqnarray*}
Observe that on  a CTD-algebra the identity is an involution.

\begin{prop}\label{invol} Let $R$ be an associative algebra equipped with an involution $\iota$, i.e.~ $\iota(ab)=\iota(b)\iota(a)$. Extending $\iota$ to $T(R)$ by
$$\iota(a_1\cdots a_n) = \iota(a_1) \cdots \iota(a_n)$$
(the order of the $a_i$'s is unchanged) gives a structure of tridendriform algebra with involution to $T^q(R)$.
\end{prop}

\begin{proo} We work by induction. We compute on one hand:
\begin{eqnarray*}
\iota(ax\g by) & = & \iota(a(x*by))\\
                      & = & \iota(a)(\iota(x * by))\\
                      & = & \iota(a)(\iota(by) * \iota(x))\\
                      & = & \iota(a)(\iota(b)\iota(y) * \iota(x))
\end{eqnarray*}
and on the other hand:
\begin{eqnarray*}
\iota(by)\d \iota(ax) & = &\iota(b)\iota(y)\d \iota(a)\iota(x)\\
                      & = & \iota(a)( \iota(b)\iota(y) * \iota(x)),
\end{eqnarray*}
Hence we have proved that $\iota(ax\g by)=\iota(by)\d \iota(ax)$. The other cases are similar.
\end{proo}

\subsection{Example} Let $\To(V)$ be the tensor algebra over $V$ without constants. Equip  $\To(V)$ with the involution
$$\iota(v_1\cdots v_n) = v_nv_{n-1} \cdots v_1$$
(the order of the $v_i$'s is reversed). So $T^q(\To(V))$ is an involutive tridendriform algebra  by Proposition \ref{invol}. 

\begin{thm}\label{thminvol} The quasi-shuffle algebra $T^q(\To(V))$ is free over $V$ (equipped with the trivial involution) in the category of involutive tridendriform algebras with involution. 
\end{thm}
\begin{proo} Let $ITD(V)$ be the the free involutive tridendriform algebra on the space $V$ equipped with the trivial involution. Since $T^q(\To(V))$ is an involutive  tridendriform algebra, the inclusion $V\to \To(V) \subset T^q(\To(V))$ induces a map $ITD(V) \to T^q(\To(V))$. Then the proof is similar to the proof of Theorem \ref{Phi}. In this variation the multilinear space $T^q\To(n)$ is spanned by the ordered-ordered partitions of $\{1,\ldots , n\}$ and so its dimension is $2^{n-1} n!$.
\end{proo}

%\thebibliography
% \bibitem{} \emph{}, 

\end{document}